\documentstyle{amsppt}
\pagewidth{6.0in}
\vsize8.0in
\parindent=6mm
\parskip=3pt
\baselineskip=16pt
\tolerance=10000
\hbadness=500
\widestnumber\key{SSS}

\define\R{\Bbb R}

\topmatter

\title
Averages over curves with torsion
\endtitle

\author
Daniel Oberlin, Hart F. Smith, and Christopher D. Sogge
\endauthor

\thanks
All authors are partially supported by the NSF.
\endthanks



\abstract
We establish $L^p$ Sobolev mapping properties for averages over certain
curves in $\R^3$, which improve upon the estimates obtained by
$L^2-L^\infty$ interpolation.
\endabstract

\endtopmatter

\document
Let $T$ be the operator given by convolution in $\R^3$ against
a smooth cutoff of arclength measure on the helix 
$\,\gamma(t)=\bigl(\cos t,\sin t, t\bigr)\,,$
$$
Tf(x)=\int f(x_1-\cos t\,,\,x_2-\sin t\,,\,x_3-t)\,\phi(t)\,dt\,.
$$
For $1<p<\infty$, let $H^{s,p}(\R^3)$ denote the nonhomogeneous
Sobolev space consisting of functions in $L^p(\R^3)$ whose
fractional derivative of order $s$ also lies in $L^p(\R^3)$.
We consider the following question:
$$
\text{\it For which values of s (depending on p) does}\;\;
T:L^p(\R^3)\rightarrow H^{s,p}(\R^3)\,\text{\it ?}
$$
By duality, it suffices to consider $2\le p<\infty$. As shown
by the first two authors in \cite{OS}, a necessary condition is that
$$
\matrix\format\r\;\;&\c&\;\;\l\\
s & \le & \frac 16 +\frac 1{3p}\quad\text{if}\quad 2\le p\le 4\,,\\
\\
s & \le & \frac 1p \quad\text{if}\quad 4\le p<\infty\,.
\endmatrix
$$
Simple arguments (see for example the lemma below) show that
$T:L^2(\R^3)\rightarrow H^{\frac 13,2}(\R^3)\,.$ Interpolation
with the trivial $L^\infty(\R^3)$ boundedness of $T$ yields
a sufficient condition of $s\le\frac 2{3p}\,.$ In
particular, interpolation yields
$$
T:L^4(\R^3)\rightarrow H^{\frac 16,4}(\R^3)\,.
\tag1$$
In this note, we combine the arguments of \cite{OS} with
Bourgain's \cite{B} improvement of the conic square function estimate
of Mockenhaupt \cite{M} to obtain the following.

\proclaim{Theorem}
There exists $\sigma>0$ such that 
$$
T:L^4(\R^3)\rightarrow H^{\frac 16+\sigma,4}(\R^3)\,.
\tag2$$
\endproclaim

We should point out that $T$ is a model for curve-averaging operators
whose canonical relations have two-sided Whitney folds.  In two
dimensions these operators are much easier to analyze and optimal
results are known.  See e.g., \cite{SS} and \cite{SW}.

In three dimensions, the translation invariant operators of this type
are the averages over curves with non-vanishing torsion (a curve
$\gamma(t)$ has non-vanishing torsion if the vectors 
$\bigl\{\gamma'(t),\gamma''(t),\gamma'''(t)\bigr\}$ are linearly
independent for each $t$.) The helix and the twisted cubic, 
$\gamma(t)=(t,t^2,t^3)\,,$ are basic examples. We restrict attention
here to the helix since this operator has the light cone in $\xi$
as its folding set. A modification of Bourgain's estimate to
conic hypersurfaces with one non-vanishing principle
curvature would yield the theorem for general curves with torsion.

The value of $\sigma$ is related to the exponent $\tau$ in equation
$(132)$ of \cite{B}, which is not explicitly determined. Any
$\sigma<\frac 13\tau$ works. 
In particular, an optimal value $\tau=\frac 14$
would yield the nearly optimal condition $\sigma<\frac 1{12}.$
Recently, Tao and Vargas \cite{TV} have modified Bourgain's arguments and
obtained a definite value of $\tau$. The authors would like to thank
T. Tao for a helpful conversation regarding Bourgain's work.

To begin the proof of $(2)$, let
$$\widehat{T}(\xi)=\int e^{-i\xi_1\cos t-i\xi_2\sin t-i\xi_3\,t}\,
\phi(t)\,dt
\tag3$$
denote the Fourier multiplier associated to $T$.

Let $\xi'=(\xi_1,\xi_2)\,.$
The oscillatory integral $(3)$ has no critical points for
$|\xi'|<|\xi_3|$. The following thus holds.
$$
\bigl|\widehat{T}(\xi)\bigr|=\Cal O(|\xi|^{-N})\quad\forall N\,,
\quad\text{if}\quad |\xi'|\le .99\,|\xi_3|\,.
$$

For $|\xi'|>|\xi_3|$ there are two, nondegenerate critical
points. The following is thus a consequence of Van der Corput's Lemma,
$$
\bigl|\widehat{T}(\xi)\bigr|\le C\,|\xi|^{-\frac 12}\,,
\quad\text{if}\quad |\xi'|\ge 1.01\, |\xi_3|\,.
$$

A simple interpolation argument implies $(2)$ for the operator
obtained by conicly restricting $\widehat T(\xi)$ to either of the
above regions. Indeed, since these bounds imply that these two
localized pieces gain a $1/2$-derivative on $L^2$, the interpolation
argument behind (1) yields estimates of the form (2) for each term
with the desired $\sigma=1/12$.

It thus suffices to establish $(2)$ for the operator
$S$ obtained by restricting the multiplier $\widehat T(\xi)$ to the
region $A$, defined by $.98\le |\xi'|\,/\,|\xi_3|\le 1.02$, 
via a smooth conic cutoff. Let $S_\lambda$ denote the operator obtained
by further restricting to the region $\lambda\le |\xi_3|\le 2\lambda\,.$
The theorem is then a result of  showing that, for some number $a>0$,
for all $\lambda>2$,
$$
\|S_\lambda\|_{4,4}\le C\,\bigl(\log\lambda\bigr)^a\,
\lambda^{-\frac 16-\frac \tau 3}\,.
\tag4
$$
We restrict attention to $\xi_3>0\,.$
Following \cite{OS}, we make a further decomposition of $S_\lambda$
by decomposing the conic set $A$ into a union of conic sets $A^j_\lambda$
as follows:
$$
\matrix\format\r&\,\c\,&\l\\
\text{for $j\ge 1$, set } A_\lambda^j & = & 
\{ 1+2^{j-1}\,\lambda^{-\frac 23} \leq |\xi'|\,/\,\xi_3 
\leq 1+2^j\,\lambda^{-\frac 23} \}\,;\\
\\
\text{set }A_\lambda^0 & = & 
\{ 1-\lambda^{-\frac 23} \leq |\xi'|\,/\,\xi_3 
\leq 1+\lambda^{-\frac 23} \}\,;\\
\\
\text{for $j\le -1$, set } A_\lambda^j & = &
\{ 1-2^{|j|}\,\lambda^{-\frac 23} \leq |\xi'|\,/\,\xi_3 
\leq 1-2^{|j|-1}\,\lambda^{-\frac 23} \}\,.
\endmatrix
$$
Introducing a suitable partition of unity on the Fourier transform
side leads to the decomposition
$$
S_\lambda= \sum_j S_\lambda^j.
$$
Inequality (4) will follow from
$$
\|S_\lambda^j\|_{4,4} \leq C\,\bigl(\log \lambda\bigr)^a
\,\lambda^{-\frac 16-\frac \tau 3}\,
2^{\frac{|j|}2(\tau-\frac 14)}
\tag5
$$
for all $j$ and $\lambda$. At this point we make a further
decomposition as in \cite{M}
of $A_\lambda^j$ into sets $A_\lambda^{jm}$ supported in $\xi'$ sectors
of angle
$\delta\doteq 2^{|j|/2}\lambda^{-\frac 13}$. This leads to a decomposition
$$
S_\lambda^j=\sum_{m=1}^{\delta^{-1}} S^{jm}_\lambda.
$$
In the notation of Theorem 1.0 of \cite{M}, we have 
$$
\widehat S^{jm}_\lambda(\xi)=
\widehat\psi_m\bigl(\lambda^{-1}\xi',\lambda^{-1}(1+\delta^2)\,\xi_3\bigr)\,
\widehat T(\xi)\,.
$$
The quantity $N$ of that theorem is related
to $j$ and $\lambda$ by $N=\delta^{-1}\,.$
\proclaim{Lemma}
$$
\|S_\lambda^{jm}\|_{4,4}\le C\,\lambda^{-\frac 14}\,\delta^{\frac 14}\,.
$$
\endproclaim
\demo{Proof}
The proof is almost identical to that of the Lemma in \cite{OS},
and is obtained by interpolating the following estimates
$$
\matrix
\format\r&\,\c\,&\l\\
\|S_\lambda^{jm}\|_{2,2} & \leq & C\,(\lambda\delta)^{-\frac 12}\,, \\
\\
\|S_\lambda^{jm}\|_{\infty,\infty} & \leq & C\,\delta\,.
\endmatrix
\tag6
$$
The first estimate in (6) is the bound 
$|\widehat S^j_\lambda(\xi)|\le C\,(\lambda\,\delta)^{-\frac 12}$, which
follows from Van der Corput's Lemma as shown in \cite{OS}. 
For the second estimate,
we consider the term $m$ corresponding to the $\xi'$ sector
along the negative $\xi_2$ axis. The convolution kernel of $S^{jm}_\lambda$,
written in the new coordinates
$$
\bigl(y_1,y_2,y_3\bigr)=\bigl(x_1,x_2+\alpha x_3,\alpha x_3-x_2\bigr)\,,
\qquad \alpha=(1+\delta^2)^{-1}\,,
$$
takes the form
$$
K(y)=\lambda^3\,\delta^3\,\int
\phi(t)\,\theta\bigl(\lambda\,\delta\,(y_1-\cos t)\,,\,
\lambda\,\delta^2\,(y_2-\sin t-\alpha t)\,,\,
\lambda\,(y_3+\sin t-\alpha t)\bigr)\,dt\,.
$$
Here and below, $\theta$ 
denotes a Schwartz function with seminorms bounded independent
of $j,m,\lambda\,,$ and with $\widehat\theta(\eta)=0$ for $\eta_3\le 1$.
We need to show that $\|K\|_{L^1}\le C\,\delta\,,$ 
and may thus replace
$\phi(t)$ by $\phi_\delta(t)$ which vanishes for $|t|\le 10\,\delta\,.$
We write $\theta=\partial_3\theta$ for some new $\theta$ to express $K(y)$
as
$$
\lambda^2\delta^3\int 
\left(\frac{\phi_\delta(t)}{\alpha-\cos t}\right)'\,\theta(\cdots)\,dt
+
\lambda^3\delta^4\int 
\frac{\sin t\,\phi_\delta(t)\,\theta(\cdots)}{\alpha-\cos t}\,dt
+
\lambda^3\delta^5\int
\frac{(\alpha+\cos t)\,\phi_\delta(t)\,\theta(\cdots)}{\alpha-\cos t}\,dt
\,.
$$
The inequality $\alpha-\cos t\ge t^2/10$ for $|t|\in[10\,\delta,\pi]$,
together with $|\phi'_\delta(t)|\le C\,\delta^{-1}\le C\,\lambda^{1/3}$,
yields the desired $L^1(dy)$ norm bounds on the first and third terms.
The desired bound for the second term follows by a further integration
by parts of the same kind.
\qed
\enddemo

To conclude the proof of $(5)$, we apply Bourgain's estimate $(132)$
of \cite{B} to obtain
$$
\Bigl\|\;\sum _m S_\lambda^{jm}f\;\Bigr\|_4 \leq
C\,\delta^{\tau-\frac 14}\,
\Bigl\|\,\Bigl(\,\sum_m |S_\lambda^{jm}f|^2\,\Bigr)^{\frac
12}\,\Bigr\|_4
\,.
$$
The number of indices $m$ is $O(\delta^{-1})$, so
$$
\sum_m\,\bigl|S_\lambda^{jm}f(x)\bigr|^2 \leq C\,\delta^{-\frac 12}\,
\Bigl(\,\sum_m \,\bigl|S_\lambda^{jm}f(x)\bigr|^4\,\Bigr)^{\frac 12}.
$$
With $\widehat f_m$ representing the localisation of ${\widehat f }$
to an appropriate sector in $\xi'$, we thus have
$$
\matrix\format\r&\quad\c\quad&\l\\
\Bigl\|\;\sum \limits _m S_\lambda^{jm}f\;\Bigr\|_4 & \leq &
C\,\delta^{\tau-\frac 12}\,
\Bigl\|\,\Bigl(\, \sum \limits _m \,|S_\lambda^{jm}f|^4 \,\Bigr)^{\frac 14}\,
\Bigr\| _4\\
\\
{} & \le & 
C\,\lambda^{-\frac 16 -\frac \tau 3}\,2^{\frac{|j|}2(\tau-\frac 14)}\,
\Bigl\|\,\Bigl(\,\sum \limits _m\, |f_m|^4 \,\Bigr)^{\frac 14}\,\Bigr\|_4\\
\\
{} & \le &
C\,\lambda^{-\frac 16 -\frac \tau 3}\,2^{\frac{|j|}2(\tau-\frac 14)}\,
\Bigl\|\,\Bigl(\, \sum \limits _m \,|f_m|^2 \,\Bigr)^{\frac 12}\,\Bigr\|_4\,.
\endmatrix
$$
A result of C\'ordoba [C] gives
$$
\Bigl\|\,\Bigl(\, \sum \limits _m\, |f_m|^2\Bigr)^{\frac 12} \Bigr\|_4
\leq C\,|\log\delta\,|^a\, \|f\|_4
$$
for some positive $a$, which completes the proof of $(5)$.

\enddocument